\documentclass[conference,a4paper]{IEEEtran}
\IEEEoverridecommandlockouts

\usepackage{graphicx}
\usepackage{cite}

\usepackage{hyperref}
\hypersetup{%
plainpages=false,
colorlinks,urlcolor=blue,linkcolor=blue,citecolor=blue,
pdftitle={Capacity Value of Additional Generation: Probability Theory
  and Sampling Uncertainty},
pdfauthor={ Chris Dent and Stan Zachary},
pdfsubject={},
pdfpagemode={UseOutlines},
pdfpagelayout={OneColumn},
pdfstartview={FitBH}
}

\newcommand{\elcc}{\nu^{\mathrm{ELCC}}}
\newcommand{\efc}{\nu^{\mathrm{EFC}}}

\usepackage{amsmath}
\usepackage{amsthm}
\usepackage{amsfonts}
\usepackage{graphicx}
\usepackage{color}
\usepackage{cite}

\renewcommand{\Pr}{\mathbf{P}}
\newcommand{\E}{\mathbf{E}}

\theoremstyle{remark}

\begin{document}

\title{Capacity Value of Additional Generation: Probability Theory and
  Sampling Uncertainty}

\author{\IEEEauthorblockN{C.J.~Dent\IEEEauthorrefmark{1},
    S.~Zachary\IEEEauthorrefmark{2}}
  \IEEEauthorblockA{\IEEEauthorrefmark{1}School of Engineering and
    Computing Sciences, Durham University, Durham, UK, \\ Email:
    chris.dent@durham.ac.uk}
  \IEEEauthorblockA{\IEEEauthorrefmark{2}School of Mathematical and
    Computer Sciences, Heriot-Watt University, Edinburgh, UK, \\
    Email: s.zachary@hw.ac.uk} \thanks{CJD was supported by the
    Supergen FlexNet consortium as part of the UK Research Councils
    Energy Programme under grant [EP/E04011X/1], and by an associated
    Engineering and Physical Sciences Research Council Knowledge
    Transfer Secondment to National Grid under grant
    [EP/H500340/1]. SZ was supported by a Scottish Funding Council
    SPIRIT knowledge-exchange grant and by an Engineering and Physical
    Sciences Research Council research grant [EP/I017054/1].} }

\maketitle

\begin{abstract}
  The concept of capacity value is widely used to quantify the
  contribution of additional generation (most notably renewables)
  within generation adequacy assessments.  
  This paper surveys the existing probability theory of assessment
of the capacity value of additional generation, and discusses the
available statistical estimation methods for risk measures which
depend on the joint distribution of demand and available additional
capacity (with particular reference to wind).

  Preliminary results are presented on assessment of sampling
  uncertainty in hindcast LOLE and capacity value calculations, using
  bootstrap resampling.  These results indicate strongly that, if
  the hindcast calculation is dominated by extremes of demand minus wind, 
  there is very large sampling
  uncertainty in the results due to very limited
  historic experience of high demands coincident with poor wind
  resource.  For meaningful calculations, some form of statistical
  smoothing will therefore be required in distribution estimation.
\end{abstract}

\begin{IEEEkeywords}
  Power system planning, Power system operation, Power system
  reliability, Risk analysis, Wind energy
\end{IEEEkeywords}

\IEEEpeerreviewmaketitle

\section{Introduction}

{T}{he} concept of capacity value is widely used to quantify the
contribution of variable output renewable generation technologies
within generation adequacy assessments.  Common specific definitions
include Effective Load Carrying Capability (ELCC, the extra demand
which the additional generation can support without increasing the
chosen risk metric), and Equivalent Firm Capacity (EFC, the completely
firm conventional capacity which would give the same risk level if it
replaced the additional variable generation).  These are usually
calculated with respect to adequacy risk indices such as the Loss of
Load Probability (LOLP) at time of annual peak, or the Loss of Load
Expectation (LOLE, the sum over periods of LOLP, or equivalently the
expected number of periods of shortage in a given time window).

There are many surveys in the literature of capacity value calculation
methods, for instance \cite{tf09,amelin,ivgtf}. In addition, a number
of papers have been published recently on analytical calculation
approaches which are valid for small additional capacities
\cite{zdjrr,zmethod}, or for the special case where the distribution
of margin of existing capacity over demand has an exponential tail
\cite{garver,windgarver}; these analytical approaches are surveyed in
\cite{dentsimp}. \cite{banda} and \cite{libook} provide general
surveys of adequacy assessment methods. The website of the IEEE PES
LOLE Working Group contains many useful presentations on current
industrial adequacy assessment practices \cite{LOLEWG}.

This paper will provide a comprehensive survey of the existing
probability theory of capacity value calculations (Section
\ref{sec:prob-theory-capac}).  In particular, this formulates the
theory in terms of the distributions of available additional capacity
and of margin of available conventional capacity over demand, which
simplifies much of the mathematical exposition and clarifies exactly
what features of distributions drive the capacity value results.
Section~\ref{sec:stat-estim} then discusses the statistical methods
for estimating the inputs to these calculations, and in particular
proposes bootstrap resampling as a means of estimating sampling
uncertainty in capacity value and LOLE calculation results. Data for
illustrative examples is provided in Section \ref{sect:data}, then
Section \ref{sect:unc} presents examples of uncertainty
assessment. While we do not claim to be presenting a quantitative
adequacy assessment for the Great Britain system, some conclusions may
be drawn regarding the degree of uncertainty in results produce by the
common hindcast approach. Finally conclusions are presented in Section
\ref{sect:conc}.


\section{Probability Theory of Capacity Values}
\label{sec:prob-theory-capac}

We present first a ``snapshot'' picture of the theory, appropriate to
the distributions of the variables involved at a given instant of
time.  We then consider how this generalises to extended periods of
time such a year. 

\subsection{Snapshot Theory: Definitions}
\label{sec:snapshot-theory}

Suppose that existing capacity less demand is represented by a random
variable $M$, with distribution function $F_M$ and density function
$f_M$.

The capacity value of additional generation represented by a random
variable $Y$ is in some appropriate sense a deterministic capacity
which is equivalent to it in terms of an associated risk.  In general
it may be viewed as the mean of $Y$ less some correction which
corresponds to its variability.

Suppose that $Y\ge0$ with distribution function $F_Y$ and density
function $f_Y$.  We denote its mean and variance by $\mu_Y$ and
$\sigma^2_Y$ respectively.

The two most commonly used definitions of the
\emph{capacity value} of $Y$ are:

\emph{Effective Load Carrying Capability (ELCC):} This is given by the
  solution $\elcc_Y$ of
\begin{equation}
  \label{eq:21}
  \Pr(M + Y   \le \elcc_Y) = \Pr(M \le 0) = F_M(0),
\end{equation}
i.e. the amount of further demand which may be added while
maintaining the same level of risk.

\emph{Equivalent Firm Capacity (EFC):} This is given by the solution
  $\efc_Y$ of
\begin{equation}
  \label{eq:26}
  \Pr(M + Y \le 0) = \Pr(M + \efc_Y \le 0) = F_M(-\efc_Y),
\end{equation}
i.e. the amount of deterministic capacity $\efc_Y$ whose addition
would result in the same level of risk as that of the addition of the
random capacity $Y$.

It is important to note that both $\elcc_Y$ and $\efc_Y$ depend on
the distributions of both $M$ and $Y$.
Note also that in the case where $Y$ is deterministic ($Y=\mu_Y$
always) we have $\efc_Y=\elcc_Y=\mu_Y$.

We consider first the case where $M$ and $Y$ are independent, and then
discuss what modifications are required to deal with the more general
case.

\subsection{The case where $M$ and $Y$ are independent}
\label{sec:case-where-m}

Assume that $M$ and $Y$ are independent.  Consider first two special
cases, both of which are analytically tractable and inform the general
case.

\subsubsection{Small additional capacity}~

The first case is where the variation in $Y$ is small in relation to
that in $M$ and was considered in \cite{zdjrr}.  
Result~1 of that paper showed that, to a good approximation,
\begin{equation}
  \label{eq:23}
  \efc_Y =\elcc_Y= \mu_Y - \frac{f_M'(0)}{2f_M(0)} \sigma_Y^2,
\end{equation}
where the error is negligible in relation to $\sigma_Y^2$ as the
latter becomes small (in relation to the variation in $M$).   We also
see here that capacity values of small \emph{independent} additions are
additive, i.e.\ if $M$, $Y_1$ and $Y_2$ are independent random
variables, and if $Y=Y_1+Y_2$, then (since $\mu_Y=\mu_{Y_1}+\mu_{Y_2}$
and $\sigma^2_Y=\sigma^2_{Y_1}+\sigma^2_{Y_2}$) it follows from
\eqref{eq:23} that $\efc_Y=\efc_{Y_1}+\efc_{Y_2}$
(and similarly for $\elcc_Y$).

\subsubsection{Exponential left tail for $M$ (`Garver approximation')}~

The second case. which forms the basis of the well-known Garver
approximation \cite{garver,windgarver}, arises when the distribution
function $F_M$ of $M$ may be treated as exponential below some level
$m_0$, i.e.\ $F_M(m)=c\exp\lambda_Mm$ for $m\le m_0$ for some
$\lambda_M>0$.  Then, since $Y\ge0$ is independent of $M$, it is
easily checked that the distribution of $M+Y$ is similarly exponential
below the level $m_0$, i.e.\ for $m\le m_0$,
\begin{align}
  \label{eq:11}
  \Pr(M+Y \le m)
  & = \int_0^\infty\mathrm{d}yf_Y(y)F_M(m-y)\nonumber\\
  & = ce^{\lambda_Mm}\int_0^\infty\mathrm{d}yf_Y(y)e^{-\lambda_My}
  \nonumber \\
  & = \Pr(M+\nu\le m),
\end{align}
where $\nu$ is the solution of
\begin{equation}
  \label{eq:6}
  \E\exp(-\lambda_M Y)
  = \int_0^\infty\mathrm{d}yf_Y(y)e^{-\lambda_My}=\exp(-\lambda_M \nu).
\end{equation}
and $\E$ denotes expectation.

It follows that, provided we can take $m_0\ge 0$ (in the case of the
EFC) and $m_0\ge \elcc_Y$ (in the case of the ELCC), the commonly used
capacity measures $\elcc_Y$ and $\efc_Y$ are each equal to $\nu$ here.
Further, in this case the capacity value depends on the distribution
of $Y$ through its Laplace transform evaluated at $\lambda_M$, so that
we may readily study which features of this distribution influence the
capacity value.  The essentials of this derivation of the Garver
approximation are not new, however the notation
of this paper makes the working much more concise than previous
versions such as \cite{windgarver}.

Finally it follows from the standard result for the Laplace transform
of the sum of independent random variables that (as in the earlier
case where the variation in $Y$ is small) 
capacity values of \emph{independent} additions are additive.

\subsubsection{The general case}~

In the more general case, when the above exponential approximation of
$F_M$ is not necessarily available, the equation~\eqref{eq:26} may be
written as
\begin{equation}
  \label{eq:7}
  \int_0^\infty \mathrm{d}y\, f_Y(y) F_M(-y) = F_M(-\efc_Y)
\end{equation}
using the standard result for the convolution of two independent
random variables (a similar expression holds for the ELCC). This may
be solved by standard numerical techniques.

\subsection{The case where $M$ and $Y$ may be dependent}
\label{sec:case-where-m-2}

If independence of $M$ and $Y$ is not assumed, then capacity values
must typically be determined by reference to \eqref{eq:21} or
\eqref{eq:26}.  We discuss this further in
Section~\ref{sec:stat-estim}.  However, when  the
variation in the distribution of the additional capacity $Y$ is
significantly less than that in the distribution of $M$ then the
equation~\eqref{eq:26} may still be approximated by \eqref{eq:23}
provided the density $f_Y(y)$ of $Y$ is replaced by its conditional
density given $M=0$ (this corresponds to the assertion that for values
of $M$ within the critical region, i.e.\ in the neighbourhood of $M=0$,
the conditional density of $Y$ does not vary significantly).  Except
for this relatively straightforward adjustment, the preceding theory
remains as before.

\subsection{Application to extended periods of time}
\label{sec:appl-extend-peri}

Typically a capacity value is required for a period of time such as a
year, which may be represented as being composed of a sequence of much
shorter periods of time, typically of hours or half-hours, to each of
which the above ``snapshot'' picture is applicable.  We index these
shorter periods by $t$, and the random variables $M_t$ and $Y_t$ are
respectively the surplus $M$ and additional capacity $Y$ during period
$t$ (their distributions typically varying with $t$).

For simplicity consider the EFC $\efc_Y$.  The
equation~\eqref{eq:26} is now replaced by
\begin{equation}
  \label{eq:10}
    \sum_t\Pr(M_t + Y_t \le 0) = \sum_t \Pr(M_t +\efc_Y\le 0)
\end{equation}
i.e.\ $\efc_Y$ is the deterministic completely firm generating
capacity if the same \emph{loss of load expectation} is to be
maintained as in the case with additional stochastic capacity $Y$.  In
the case where the pairs $(M_t,Y_t)$ are viewed as corresponding to
randomly chosen periods of time and are thus treated as identically
distributed (with the `whole peak season' distribution), the
equation~\eqref{eq:10} does of course reduce to \eqref{eq:26}.

\section{Statistical estimation}
\label{sec:stat-estim}

Throughout this section we consider, for definiteness, the estimation
of the EFC $\efc_Y$.

Approaches to statistical estimation naturally depend on the nature of
the available data.  In the present paper we assume that the random
variable $M=X-D$ where $X$ is existing, typically conventional,
capacity and $D$ is demand.  We further assume that the successive
instances of $X$, i.e.\ the random variables~$X_t$, may be modelled as
identically distributed, with a \emph{known} distribution function~$F_X$
(typically determined from a capacity outage table),
and are independent of the random variables $D_t$ and $Y_t$.  We thus
rewrite the equation~\eqref{eq:10} as
\begin{equation}
  \label{eq:1}
  \sum_t \E F_X(D_t - Y_t) = \sum_t \E F_X(D_t -\efc_Y).
\end{equation}
This formulation incorporates the above independence assumption, and
is thus most suited to inference given observations
$(d_t,y_t)$ of the successive pairs $(D_t,Y_t)$.  (We here assume that
the joint distribution of the pairs $(D_t,Y_t)$ is to be estimated
from data, as will be the case when the additional capacity $Y$
corresponds to some renewable resource such as wind generation.)

Further assumptions (often not made explicit in the existing
literature) are now necessary in order to make any inference for
$\efc_Y$.  Usually we treat the pairs $(D_t,\,Y_t)$ as identically
distributed---with the `whole peak season' distribution referred to in
the previous section---and we henceforth assume this, representing
this common distribution as that of a generic pair $(D,\,Y)$.  Point
estimates of $\efc_Y$ do not then in general require any assumptions
about the nature of any dependence between the successive pairs
$(D_t,Y_t)$.  Assessments of uncertainty for these estimates, e.g.\
confidence intervals, do require such assumptions.  The simplest such
is to take the pairs $(D_t,Y_t)$ to be additionally independent, thus
allowing straightforward techniques to be used for the construction of
confidence intervals, etc.  However, in reality significant serial
correlations will be present; correctly allowing for these (which may
require data of better quality than is currently available) will
significantly increase the reported uncertainty associated with
estimates of $\efc_Y$.

Finally, while it will be typically be unrealistic to assume
independence of the generic $D$ and $Y$ (since, for example, demand
and available wind will typically have some statistical association),
we make nevertheless wish to make some smoothness assumptions
concerning the nature of the dependence between them.  We discuss this
further below.

The simplest approach to the estimation of $\efc_Y$ is to substitute
the observed values $(d_t,\,y_t)$ of $(D_t,\,Y_t)$ into \eqref{eq:1},
giving the expression
\begin{equation}\label{eq:2}
  \sum_t F_X(d_t - y_t) = \sum_t F_X(d_t -\efc_Y).
\end{equation}
(where the sum is over the historic times $t$ for which data is
available) and to solve for $\efc_Y$. This is the commonly used
\emph{hindcast} approach, and has the virtue of making no assumptions
about the common joint distribution of $(D,\,Y)$.  (In particular it
attempts to account for statistical association between $D$ and $Y$ in
the historic data without requiring any advanced statistical
technology).  A common criticism of the approach relates to the
difficult of making an uncertainty assessment of the estimate for
$\efc_Y$; indeed results are usually presented without assessment of
uncertainty.  However, under the additional independent assumption for
the successive pairs $(D_t,\,Y_t)$ referred to above, the latter may
be obtained easily by bootstrapping (successively resampling from the
data) \cite{efronas,efronbook}.

A major difficulty with the simple hindcast approach described above
is the likely shortage of data in the extreme regions of relevance to
the solution of \eqref{eq:1} (or \eqref{eq:2}), leading to
considerable uncertainty in the estimate of $\efc_Y$.  This can be
partially remedied by assuming some reasonably smooth form of
dependence in the joint distribution of $(D,\,Y)$---an assumption
which effectively allows more observations to make a helpful
contribution to the estimation procedure.  The most radical such
assumption is to take $D$ and $Y$ to be independent, but this will
typically be unrealistic: for example, to the extent that wind and
demand are typically both higher in winter, there is a positive
statistical association between them.

A sensible compromise is to indeed proceed as if the demand $D$ and
additional capacity $Y$ were independent, but to replace the
distribution of $Y$ by an estimate of its \emph{conditional}
distribution given that the demand $D$ lies within the region critical
for the estimation of $\efc_Y$.  In practice this means treating $D$
and $Y$ as independent but replacing the distribution of $Y$ by its
distribution for those \emph{times}---of the year and of the day---in
which the modelled system is at risk, i.e.\ in which the demand $D$
tends to be high.  (In the case where $Y$ is wind generation, an
interesting alternative is to use a distribution of $Y$ estimated
using data from those classes of weather system for which demand~$D$
is likely to be high \cite{zdbjrr}; however care is required here as
it may be that
types of weather for which demand is less extreme may also be
associated with particularly low levels of wind, again placing the
system at risk.)  From \eqref{eq:1} and our ``identically
distributed'' assumption, estimation of $\efc_Y$ now reduces to the
solution of
\begin{equation}
  \label{eq:3}
    \E F_X(D - Y) = \E F_X(D - \efc_Y).
\end{equation}
where $D$ and $Y$ are additionally treated as independent with
distributions as suggested by the data sampled at the critical periods
referred to above.  These may be the relevant empirical distributions,
perhaps smoothed, in which case the solution of~\eqref{eq:3} will
necessarily be numerical; alternatively parametric estimates of these
distributions may be used.  In the former case assessments of
uncertainty may be made by bootstrap resampling using each of the two
empirical distributions; in the latter either bootstrap or analytical
techniques may be used.



\section{Data for Examples}\label{sect:data}

This section describes the Great Britain-based test data used for
calculation examples in this paper.  The descriptions are quite brief,
as this paper does not claim to perform a quantitative capacity value
study for this system. Instead, it illustrates the ideas on
uncertainty assessment by bootstrapping
described above, and provides an
indication of the degree of uncertainty which may be observed in more
definitive calculations.

\subsection{Conventional Plant}

The probability distribution of available conventional capacity is
based on the list of units connected to the GB system in winter
2008/09\footnote{For this initial indicative study of uncertainty,
  demand will be rescaled to give a risk level which is sustainable in
  the long run. For a quantitative GB adequacy study, the conventional
  plant data would have to be updated and an appropriate projection of
  underlying demand patterns used.}.
A
standard capacity outage probability table (COPT) calculation is
performed, with the availability probabilities for each class of
generating unit taken from \cite{ngwo1011}; in a small number of
cases, the maximum contribution from each station is capped due to
finite network capacity or emissions constraints, following the
practice of the GB system operator.  In all examples, the
distributions of available conventional capacity at different times
are assumed to be identically distributed. This is reasonable if there
is little planned maintenance at times when risk is high, and hence
all units that are mechanically available are available to generate if
required.

For the illustrative examples presented here, the distribution of
available conventional capacity will not be rescaled; instead the peak
demand will be adjusted to achieve a risk level consistent with
historic experience. This will allow general conclusions to be drawn
regarding the the importance of uncertainty analysis, but a more
careful treatment of the conventional plant is clearly required in
practical risk calculations (and associated assessment of
uncertainty).


\subsection{Demand Data}

Half-hourly historic transmission-metered demand data is available for
the GB system since April 2001 \cite{demanddata}.  This paper uses the
seven years 2002-8 for which coincident wind resource data is also
available; as the wind data is hourly, for each hour the demand used
is the higher demand from the two half hours contained therein.
Historic
demands from different years may be compared to an extent by rescaling
according to each winter's \emph{Average Cold Spell} (ACS) peak demand
metric\footnote{ACS peak demand is the standard measure of underlying
  peak demand level in Great Britain, independent of the weather
  conditions in the year in question.  See the glossary of
  \cite{sqss} for a formal definition.}, although this does not
account fully for changes in
underlying demand patterns (including increasing penetrations of
distributed generation). The ACS peak level for each winter is
published by the GB System Operator.
  


\subsection{Wind Data}

The wind resource data used in this paper was generated by P\"oyry
Consulting for their `Impact of Intermittency' report; a more detailed
description of the dataset may be found in \cite{poyrypublic}. This
dataset is based on hourly wind speed records from 19 onshore meteorological
stations around GB, plus 7 offshore locations where the historic data
is derived by atmospheric modelling. As described above, these observations
are coincident in time with the demand data from 2002-8.
These wind speed time series are
transformed to wind load factors using a generic wind farm power
curve, and the resulting load factor time series from each of these
sites is assumed to be representative of wind farms in that area.
Finally a time series of system aggregate wind power outputs are
generated for a 2020 scenario of installed wind capacities.

We note a number of uncertainties associated with this data, including
how representative the chosen locations are of actual wind farm
locations, some issues of limited geographical coverage (notably
including a total lack of data from Scottish offshore waters) and
conversion of wind speed data to hub height wind speeds and then wind
power.  For quantitative applied studies, we believe that the most
satisfactory approach to wind resource modelling is mesoscale
reanalysis, in which physical atmospheric modelling is used to
downscale a coarse-grid historic record to the spatial and temporal
resolution required for the power system analysis work
\cite{lewpesgm}.

\section{Assessment of Uncertainty}\label{sect:unc}

This section demonstrates how bootstrap resampling, mentioned in
Section~\ref{sec:stat-estim}, can provide an estimate of statistical
uncertainty in hindcast LOLE and EFC calculations.  Bootstrap brings
the key advantage over the previous method in \cite{haschereq} that it
does not require an assumed `correct' result for comparison in order
to make a quantification of statistical uncertainty.  As described in
\ref{sect:data}, we do not claim that this paper provides a
quantitative LOLE or capacity value assessment for the GB
system. However, the consideration of uncertainty in our example
calculations will permit important conclusions about statistical
methodologies for practical adequacy assessment.

For this illustrative purpose, all calculations are performed with ACS
peak demand of 61.5 GW and the conventional plant distribution
described above. This underlying demand level gives, for all wind
capacities considered, a highest hourly LOLP in the hindcast LOLE calculation
(i.e. largest term in the sum on the right hand side of (\ref{eq:2})) of
the order of 1\%; a risk level of this order is widely reckoned to be
economically sustainable in GB. The
hindcast calculation utilises only those observations $(d_t,\,y_t)$
corresponding to the 5000 highest values of $d_t$, and constituting
8\% of the total data.  (For the remaining observations $d_t$ is
insufficiently high to make any contribution to the solution of
\eqref{eq:2}.)

\subsection{Uncertainty in Hindcast LOLE and EFC}

\begin{figure}%
\begin{center}
\includegraphics[width=7cm]{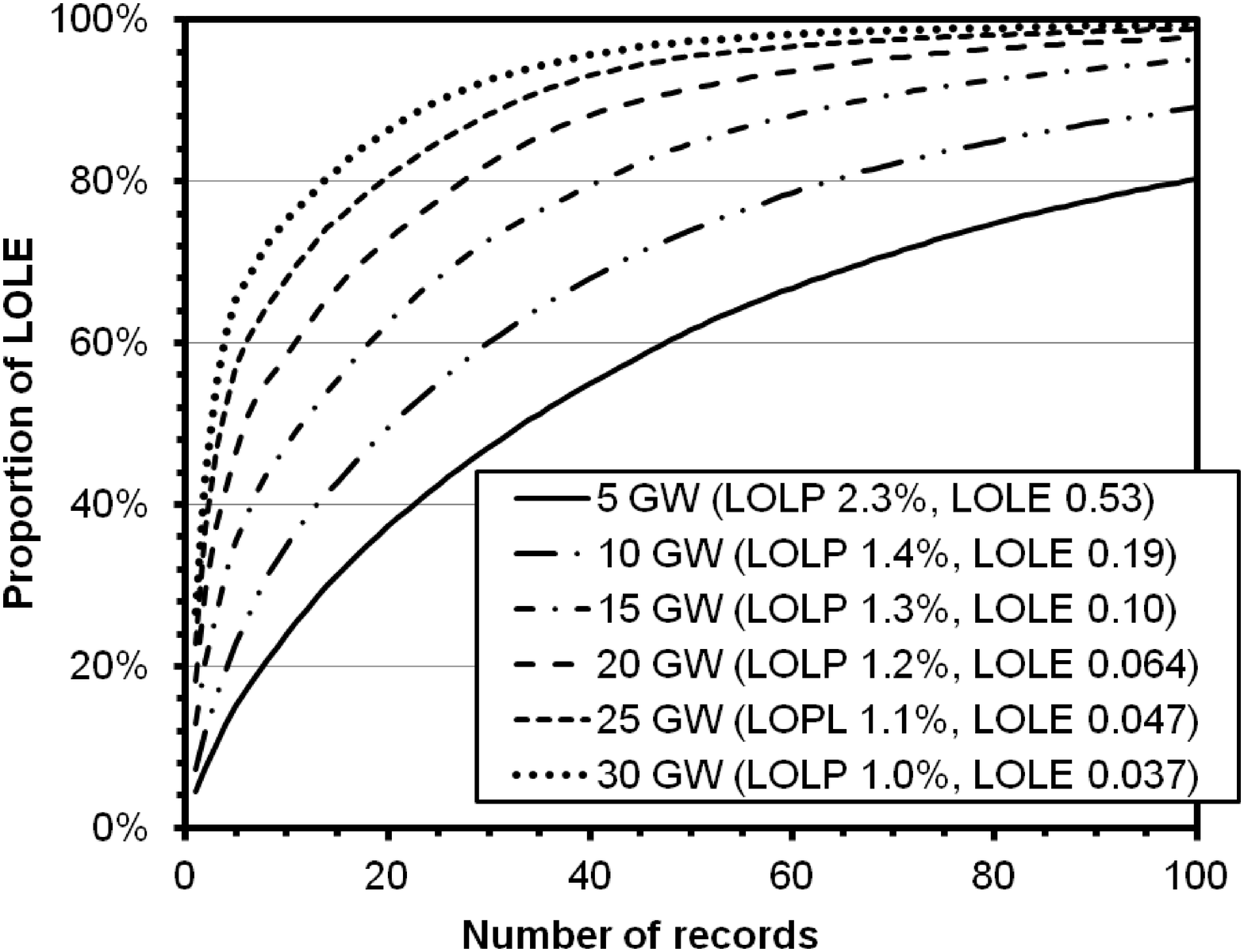}%
\end{center}
\caption{Proportion of LOLE arising from the top $n$ net demands in the right hand side of (\ref{eq:2}), for
  a range of installed wind capacities. The LOLE figures quoted are in hours per seven years (based on the seven years 
  of data used. The LOLP figure is the highest hourly LOLP in the sum on the RHS of (\ref{eq:2}).}%
\label{fig:LOLEvN}%
\end{figure}

The part of the calculated LOLE due to the $n$ historic records
with the highest net demands
is shown in Fig.~\ref{fig:LOLEvN} for a range of installed wind
capacties (i.e.
\begin{equation}
  \sum_{t=1}^n F_X (d_t - y_t)
\end{equation}
where the times $t$ are ordered by decreasing net demand $d_t-y_t$). 
This is expressed as a proportion of the total LOLE
where all data is included in the sum.

An indication of uncertainty in results for a range of wind
capacities is obtained by simple rescaling of the wind power data from the
2020 scenario on which the time series is based.
It may be seen that at high wind capacities, the calculation becomes
dominated by a very small number of records with high demand and poor
wind resource. Indeed, for 30 GW installed wind capacity, 67\% of the
calculated LOLE is due to records from just two consecutive days in
February 2006. It is clear that in such a situation the calculation
will have very limited ability to estimate risk at future times.

\begin{figure}%
\begin{center}
\includegraphics[width=7cm]{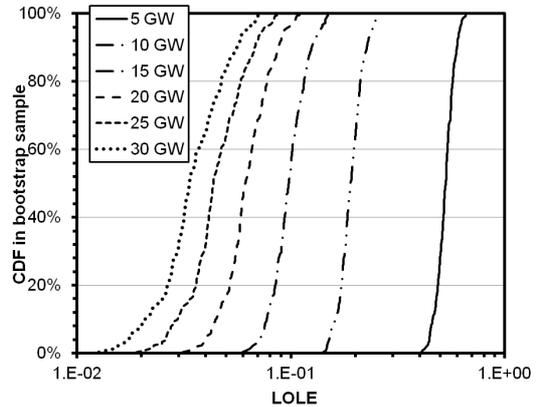}%
\end{center}
\caption{Distribution of LOLE results arising from 200 bootstrap
  samples from the historic time series of wind and demand, for a
  range of installed wind capacities.}%
\label{fig:BSlole}%
\end{figure}

Fig.~\ref{fig:BSlole} shows the distribution of calculated LOLEs
arising from 200 bootstrap samples from the historic joint series of
demand and available wind capacity.  Each bootstrap sample is a random
sample of the same size as that of the original dataset (of size 5000
as described above), sampled from the empirical distribution of the
that dataset.  The distribution of the calculated LOLEs of the
200 bootstrap samples constitutes a direct assessment of the sampling
distribution of the LOLE for the original data.  The calculations are
repeated for a range of installed wind capacities.
In making this assessments, the demand-wind pairs at different times
are assumed to be independent and identically distributed. In reality
there are important serial associations between consecutive hours and
days, and this estimate of uncertainty based on the assumption of
independence will thus be an underestimate.

Even for only 5 GW of installed wind capacity, the 95\% confidence
intervals for the LOLE arising from the bootstrap analysis ranges from
0.44 to 0.62, a factor of 1.4.  However, for 30 GW installed capacity,
the confidence interval covers the range 0.016 to 0.068, a factor of
4.2.

\begin{figure}%
\begin{center}
\includegraphics[width=7cm]{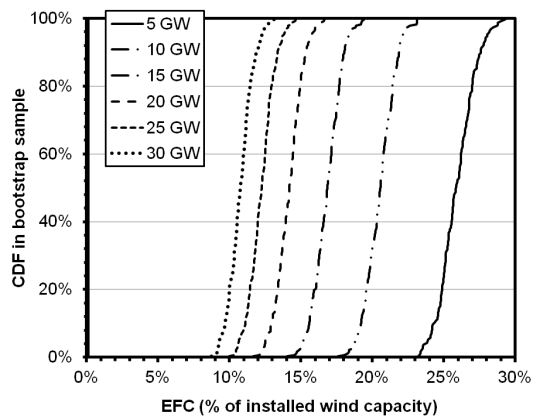}%
\end{center}
\caption{Distribution of EFC results arising from 200 bootstrap
  samples from the historic time series of wind and demand, for a
  range of installed wind capacities.}%
\label{fig:BSefc}%
\end{figure}

The distributions of EFC results arising from similar bootstrap
resampling are shown in Fig.~\ref{fig:BSefc}. 
The apparent greater robustness of the EFC results (as compared to
those for LOLE) arising from these bootstrap calculations is 
merely a consequence of the nonlinear rescaling of the risk level to a
MW capacity value by inverting the cumulative distribution function of
$M$. 

\subsection{Discussion}

As stated previously, these calculations are primarily intended to
illustrate the proposed statistical methodologies, and certainly do
not constitute a quantitative adequacy assessment for the Great Britain
system. The latter is planned by the authors as future work, and will
require improved data, including realistic scenarios of installed
generating capacity for the future years under study.

Nevertheless, some important conclusions can be drawn from this
initial study. The first is that a hindcast LOLE calculation for the
GB system with a high installed wind capacity can say very little
meaningful about real risk levels. This is demonstrated clearly by the
tiny number of historic days' data which dominate the calculation, and
by the bootstrap results. The latter in particular give an optimistic
view of the degree of uncertainty due to the modelling of demand-wind
pairs at different times as independent and identically distributed.

Section~\ref{sec:stat-estim} describes how sampling uncertainty might
be reduced by some form of statistical smoothing. This might be
relatively straightforward as described there,
or might involve application of extreme value statistical methods.
While such approaches will certainly reduce sampling uncertainty,
there will remain some (hopefully modest) additional uncertainty over
the applicability of the chosen methodology.

Future work on bootstrap assessment of uncertainty in hindcast
calculations should include consideration of serial associations in the
time series of demand and wind conditions.  This will be required for
any hindcast calculations performed in GB for relatively modest wind
capacities (where the sampling uncertainty is not overwhelming), and
also for other systems where the risk calculation results are not
dominated to the same extent by a small number of extremes of net
demand in the historic data.

\section{Conclusion}\label{sect:conc}

This paper has surveyed the existing probability theory of assessment
of the capacity value of additional generation, and discussed the
available statistical estimation methods for risk measures which
depend on the joint distribution of demand and available additional
capacity (with particular reference to wind).

Preliminary results have been presented on assessment of sampling
uncertainty in hindcast LOLE and capacity value calculations, using
bootstrap resampling.  The test system used is not sufficiently
realistic as to provide a quantitative generation adequacy risk study
for the Great Britain system.  Nevertheless, the results indicate
strongly that, in systems such as GB where the calculated adequacy
risk is dominated by extremes of demand minus wind, there is very
large sampling uncertainty in hindcast adequacy calculation results
due to limited historic experience of high demands coincident with
poor wind resource.  For meaningful calculations, some form of
statistical smoothing will therefore be required in 
estimation.  The results also confirm that uncertainty analysis is a
vital part of any generation adequacy study, particular in systems
with a substantial wind generation capacity.

\section*{Acknowledgemnts}

The authors are grateful to the National Grid Company for many
valuable discussions and for hosting Dr. Dent on a secondment where
this work was originally conceived, and to P\"oyry for the use of their
wind resource dataset.  They also acknowledge discussions
with S.\ Mokkas and K.\ Royal at Ofgem, C.\ Gibson, S.\ Mancey, members of
the IEEE PES Task Forces on Capacity Value of Wind and Solar Power,
the NERC LOLE Working Group, and colleagues at Durham and Heriot-Watt
Universities, NREL and University College Dublin.



\bibliography{jrr}{}

\begin{thebibliography}{10}
\providecommand{\url}[1]{#1}
\csname url@samestyle\endcsname
\providecommand{\newblock}{\relax}
\providecommand{\bibinfo}[2]{#2}
\providecommand{\BIBentrySTDinterwordspacing}{\spaceskip=0pt\relax}
\providecommand{\BIBentryALTinterwordstretchfactor}{4}
\providecommand{\BIBentryALTinterwordspacing}{\spaceskip=\fontdimen2\font plus
\BIBentryALTinterwordstretchfactor\fontdimen3\font minus
  \fontdimen4\font\relax}
\providecommand{\BIBforeignlanguage}[2]{{%
\expandafter\ifx\csname l@#1\endcsname\relax
\typeout{** WARNING: IEEEtran.bst: No hyphenation pattern has been}%
\typeout{** loaded for the language `#1'. Using the pattern for}%
\typeout{** the default language instead.}%
\else
\language=\csname l@#1\endcsname
\fi
#2}}
\providecommand{\BIBdecl}{\relax}
\BIBdecl

\bibitem{tf09}
A.~Keane, M.~Milligan, C.~J. Dent, B.~Hasche, {C. D'Annunzio}, K.~Dragoon,
  H.~Holttinen, N.~Samaan, L.~S\"{o}der, and {M. O'Malley}, ``{Capacity Value
  of Wind Power},'' \emph{{IEEE Trans. Power Syst.}}, vol.~26, no.~2, pp.
  564--572, 2011, {IEEE PES Task Force report.}

\bibitem{amelin}
M.~Amelin, ``{Comparison of Capacity Credit Calculation Methods for
  Conventional Power Plants and Wind Power},'' \emph{IEEE Trans. Power Syst.},
  vol.~24, no.~2, pp. 685--691, May 2009.

\bibitem{ivgtf}
{NERC Integration of Variable Generation Task Force}, ``{Methods to Model
  Calculate Capacity Contributions of Variable Generation for Resource Adequacy
  Planning},'' {2011}, available at
  \\http://www.nerc.com/filez/ivgtf\_Probabilistic.html.

\bibitem{zdjrr}
S.~Zachary and C.~J. Dent, ``{Probability theory of capacity value of
  additional generation},'' 2011, {J. Risk and Reliability (Proc. IMechE pt.
  O), in press.}

\bibitem{zmethod}
K.~Dragoon and V.~Dvortsov, ``{Z-Method for Power System Resource Adequacy
  Applications},'' \emph{IEEE Trans. Power Syst.}, vol.~21, no.~2, May 2006.

\bibitem{garver}
L.~L. Garver, ``Effective load carrying capability of generating units,''
  \emph{IEEE Trans. Power Apparatus and Systems}, vol.~85, no.~8, pp. 910--919,
  August 1966.

\bibitem{windgarver}
C.~{D'Annunzio} and S.~Santoso, ``Noniterative method to approximate the
  effective load carrying capability of a wind plant,'' \emph{IEEE Trans.
  Energy Conv.}, vol.~23, no.~2, pp. 544--550, June 2008.

\bibitem{dentsimp}
C.~J. Dent, A.~Keane, and J.~W. Bialek, ``{Simplified Methods for Renewable
  Generation Capacity Credit Calculation: A Critical Review},'' in \emph{IEEE
  PES General Meeting}, 2010.

\bibitem{banda}
R.~Billinton and R.~N. Allan, \emph{Reliability evaluation of power systems,
  2nd edition}.\hskip 1em plus 0.5em minus 0.4em\relax Plenum, 1994.

\bibitem{libook}
W.~Li, \emph{Risk Assessment of Power Systems: Models, Methods, and
  Applications}.\hskip 1em plus 0.5em minus 0.4em\relax IEEE/Wiley, 2005.

\bibitem{LOLEWG}
IEEE PES LOLE WG meeting, 7-8 November 2011, available at \\
  http://ewh.ieee.org/cmte/pes/rrpa/RRPA\_20111107LBPmtg.html.

\bibitem{efronas}
B.~Efron, ``{Bootstrap Methods: Another Look at the Jackknife},'' \emph{Annals
  of Statistics}, vol.~7, no.~1, pp. 1--26, 1979.

\bibitem{efronbook}
B.~Efron and R.~Tibshirani, \emph{An introduction to the bootstrap}.\hskip 1em
  plus 0.5em minus 0.4em\relax Chapman and Hall, 1993.

\bibitem{zdbjrr}
D.~J. Brayshaw, C.~J. Dent, and S.~Zachary, ``{Wind generation's contribution
  to supporting peak electricity demand: meteorological insights},'' 2011, {J.
  Risk and Reliability (Proc. IMechE pt. O), in press.}

\bibitem{ngwo1011}
{National Grid}, ``{Winter Outlook Report 2010/11},'' 7 October 2010, available
  from http://www.nationalgrid.com/uk/Gas/TYS/outlook/.

\bibitem{demanddata}
``{National Grid: Demand Data},'' available at
  \\http://www.nationalgrid.com/uk/Electricity/Data/Demand+Data/.

\bibitem{sqss}
``{Great Britain National Electricity Transmission System Security and Quality
  of Supply Standard, Issue 2.0},'' 24 June 2009, available from
  http://www.nationalgrid.com/uk/Electricity/Codes/\\gbsqsscode/DocLibrary/.

\bibitem{poyrypublic}
{P\"{o}yry Energy Consulting}, ``{Impact of Intermittency: Summary Report},''
  July 2009.

\bibitem{lewpesgm}
D.~Lew, C.~Alonge, M.~Brower, J.~Frank, L.~Freeman, K.~Orwig, C.~Potter, and
  {Y.-H. Wan}, ``{Wind Data Inputs for Regional Integration Studies},'' in
  \emph{IEEE PES General Meeting}, 2011, p. 1097.

\bibitem{haschereq}
B.~Hasche, A.~Keane, and M.~O'Malley, ``{Capacity Value of Wind Power,
  Calculation and Data Requirements: the Irish Power System Case},'' \emph{IEEE
  Trans. Power Syst.}, in press.

\end{thebibliography}
\bibliographystyle{IEEEtran}

\end{document}